\documentclass[12pt]{amsart}
\usepackage{a4wide, amssymb, bbm, mathrsfs, stmaryrd}
\usepackage[2cell,arrow,curve,matrix]{xy}
\UseTwocells

\title{Twisted Hochschild homology and MacLane homology}

\author{Teimuraz Pirashvili}\thanks{$^{*}$ Partially supported by
the grant from  MEC,  MTM2006-15338-C02-02 (European FEDER support
included)}
\theoremstyle{definition}

\let\Gg\Gamma

\def\Z{{\mathbb Z}}
\def\N{{\mathbb N}}

\let\onto\twoheadrightarrow

\def\cok{\mathop{\sf Coker}\nolimits}

\def \ker{\mathop{\sf Ker}\nolimits}
\def \hom{\mathop{\sf Hom}\nolimits}
\def \ext{\mathop{\sf Ext}\nolimits}

\def \tor{\mathop{\sf Tor}\nolimits}
\def\H{{\mathsf{H}}}%
\def\HML{{\mathsf{HML}}}%
\newtheorem{Pro}{Proposition}
\newtheorem{Le}[Pro]{Lemma}
\newtheorem{Th}[Pro]{Theorem}
\newtheorem{Co}[Pro]{Corollary}

\newtheorem{Exm}[Pro]{Example}

\def\sF{{\mathscr F}}
\def\sP{{\mathscr P}}

\let\bu\bullet

\def\F{{\mathbb F}}

\let\t\otimes

\begin{document}

%\begin{abstract}    % type your abstract below

%\end{abstract}
\maketitle

\begin{abstract}
We prove that $\H_i(A,\Phi(A))=0$, $i>0$. Here $A$ is a commutative
algebra over the prime field $\F_p$ of characteristic $p>0$ and
$\Phi(A)$ is $A$ considered as a bimodule, where the left
multiplication is the usual one, while the right multiplication is
given via Frobenius endomorphism and $\H_\bu$ denotes the Hochschild
homology over $\F_p$. This result has implications in MacLane
homology theory. Among other results, we prove that
$\HML_\bu(A,T)=0$, provided $A$ is an algebra over a field $K$ of
characteristic $p>0$ and $T$ is a strict homogeneous polynomial
functor of degree $d$ with $1<d<{\sf Card}(K)$.
\end{abstract}

%\maketitle

\section{Introduction}
In this short note we study Hochschild and MacLane homology of
commutative algebras over the prime field $\F_p$ of characteristic
$p>0$. Let us recall that MacLane homology is isomorphic to
 the topological Hochschild homology \cite{PW} and to
the stable $K$-theory as well \cite{4}.

 Let $A$ be a commutative algebra over the prime field $\F_p$ of
characteristic $p>0$ and let $\Phi(A)$ be denote $A$-$A$-bimodule,
which is $A$ as a left $A$-module, while the right multiplication is
given via Frobenius endomorphism. We prove that the Hochschild
homology vanishes $\H_i(A,\Phi(A))=0$, $i>0$. The proof makes use a
simple result on homotopy groups of simplicial rings, which says
that if $R_\bu$ is a simplicial ring such that all rings involved in
$R_\bu$ satisfy $x^m=x$, $m\geq 2$ identity then $\pi_i(R_\bu)=0$
for all $i>0$.  These results has implications in MacLane homology
theory. We extend the computation of Franjou-Lannes and Schwartz
\cite{FLS} of MacLane (co)homology of finite fields with
coefficients in symmetric $S^d$ and divided powers $\Gamma^d$ to
arbitrary commutative $\F_p$-algebras, provided that $d>1$. As a
consequence of our computations we show that $\HML_\bu(A,T)=0$,
provided $T$ is a strict homogeneous polynomial functor of degree
$d>1$ and $A$ is an algebra over a  field $K$ of characteristic
$p>0$ with ${\sf Card}(K)>d$.

Many thanks to referee for his valuable comments and to Petter
Andreas Bergh for his remarks.

\section{When it is too easy to compute homotopy groups} It is
well-known that the homotopy groups of a simplicial abelian group
$(A_\bu,\partial_\bu,s_\bu)$ can be computed as the homology of the
normalized chain  complex $(N_\bu(A_\bu),d)$, where
$$N_n(A_\bu)=\{x\in A_n| \partial_i(x)=0, i>0\}$$
and the  boundary map $N_n(A_\bu)\to N_{n-1}(A_\bu)$ is induced by
$\partial_0$. Our first result shows that if $A_\bu$ has a
simplicial ring structure and the rings involved in $A_\bu$ satisfy
extra conditions then homotopy groups are zero in positive
dimensions. This fact is an easy consequence of the following result
which is probably well-known

\begin{Le} Let $R_\bu$ be a simplicial object in the category of not
necessarily associative rings and let $x,y\in N_n(R_\bu)$ be two
elements. Assume $n>0$ and  $x$ is a cycle. Then the cycle $xy\in
N_{n}(R_\bu)$  is a coboundary.
\end{Le}

\begin{proof} Consider the element
$$z=s_0(xy)-s_1(x)s_0(y).$$
Then we have
$$\partial_0(z)=xy-(s_0\partial_0(x))y=xy.$$ Moreover,
$$\partial_1(z)= xy-xy=0,$$
We also have
$$\partial_2(z)= (s_0\partial_1(x))(s_0\partial_1(y))-x(s_0\partial_1(y))=0.$$
Similarly for all $i>2$ we have
$$\partial_i(z)=(s_0\partial_{i-1}(x))(s_0\partial_{i-1}(y))-(s_1\partial_{i-1}(x))(s_0\partial_{i-1}(y))=0.$$
Hence $z$ is an element of $N_{n+1}(R_\bu)$ with $\partial(z)=xy$.
\end{proof}

\begin{Co}\label{homonulia} Let  $R_\bu$ be a simplicial ring. If
the rings involved in $R_\bu$ satisfy $x^m=x$ identity for $m\geq
2$, then
$$\pi_n(R_\bu)=0, \ \ n>0.$$
\end{Co}

\begin{proof} Take a cycle $x\in N_n(R_\bu)$, $n>0$. Then the class of $x=xx^{m-1}$
in $\pi_n(R_\bu)$ is zero.
\end{proof}

{\bf Remark}. A more general fact is true. Let $\bf T$ be a pointed
algebraic theory \cite{schwede} and let $X_\bu$ be a simplicial
object  in the category of ${\bf T}$-models \cite{schwede}. Then
$\pi_1(X_\bu)$  is a group object in the category of ${\bf
T}$-models, while $\pi_i(X_\bu)$ are abelian group objects in the
category of ${\bf T}$-models for all $i>1$. Thus $\pi_i(X_\bu)=0$,
$i\geq 1$ provided  all group objects are trivial. This is what
happens for the category of rings satisfying the identity $x^m=x$,
$m\geq 2$. Another interesting case is the category of Heyting
algebras \cite{leo}.

\section{Hochschild homology with twisted coefficients}
In what follows the ground field is the prime field $\F_p$ of
characteristic $p>0$. All algebras are taken over $\F_p$ and they
are assumed to be associative. For an algebra $R$ and an
$R$-$R$-bimodule $B$ we let $\H_\bu(R,B)$ and $\H^\bu(R,B)$ be the
Hochschild homology and cohomology of $R$ with coefficients in $B$.
Let us recall that
$$\H_\bu(R,B)=\tor^{R\t R^{op}}_\bu(R,B)$$
and $$\H^\bu(R,B)=\ext_{R\t R^{op}}^\bu(R,B).$$ Moreover, we let
$C_\bu(R,B)$  be the standard simplicial vector space computing
Hochschild homology
$$\pi_\bu(C_\bu(R,B))\cong \H_\bu(R,B).$$
Recall that $C_n(R,B)=B\t R^{\t n}$, while
$$\partial_0(b,r_1,\cdots, r_n)=(br_1,\cdots, r_n),$$
$$\partial_i(b,r_1,\cdots,
r_n)=(b,r_1,\cdots, r_ir_{i+1},\cdots, r_n), \ \ 0<i<n$$ and
$$\partial_n(b, r_1,\cdots,
r_n)=(r_nb,r_1,\cdots,  r_{n-1}).$$ Here $b\in B$ and
$r_1,\cdots,r_n\in R$.

Let $n>1$ be a natural number and  let $A$ be a commutative
$\F_p$-algebra. The Frobenius homomorphism gives rise to the
functors $\Phi^n$ from the category of $A$-modules to the category
of $A$-$A$-bimodules, which are defined as follows. For an
$A$-module $M$ the bimodule $\Phi^n(M)$ coincides with $M$ as a left
$A$-module, while the right $A$-module structure on $\Phi^n(M)$ is
given by
$$ma=a^{p^n}m, \ \ a\in A,\, m\in M.$$
Having $A$-$A$-bimodule $\Phi^n(M)$ we can consider the Hochschild
homology $\H_\bu(A,\Phi^n(M))$. In this section we study these
homologies. In order to state our results we need some notation. We
let $\psi^n(A)$ be the quotient ring $A/(a-a^{p^n})$, $n\geq 1$
which is considered as an $A$-module via the quotient map $A\onto
\psi^n(A)$. Thus $\psi^n$ is the left adjoint of the inclusion of
the category of commutative $\F_p$-algebras with identity $x^m=x,
m=p^n$ to the category of all commutative $\F_p$-algebras.

\begin{Exm} Let $n>1$. If $K$ is a finite field with $q=p^d$ element then
$\psi^n(K)=K$ if $n=dt$, $t\in \N$ and $\psi^n(K)=0$ if $n\not =dt$,
$t\in \N$.
\end{Exm}

\begin{Le}\label{us} Let $A$ is a commutative algebra over a  field
$K$ of characteristic $p>0$ with ${\sf Card}(K)>p^n$. Then
$\psi^n(A)=0$, $n>1$.
\end{Le}

\begin{proof} By assumption  there exists $k\in K$ such that
$k^{p^n}-k$ is an invertible element of $K$. It follows then that
the elements of the form  $a^{p^n}-a$ generates whole $A$.
\end{proof}

\begin{Th}\label{respsi} Let $A$ be a commutative $\F_p$-algebra and $n>1$.  Then
$$\H_i(A,\Phi^n(A))=0$$ for all $i>0$ and
$$\H_0(A,\Phi^n(A))\cong \psi^n(A).$$
\end{Th}

\begin{proof} The proof consists of three  steps.

{\bf Step 1. The theorem holds if  $A=\F_p[x]$}. In this case we
have the following projective resolution of $A$ over $A\t
A=\F_p[x,y]$:
$$0\to \F_p[x,y] \buildrel \eta \over \to \F_p[x,y] \buildrel \epsilon \over
\to \F_p[x]\to 0.$$ Here $\epsilon(x)=\epsilon(y)=x$ and $\eta$ is
induced by  multiplication by $(x-y)$. Hence for any
$A$-$A$-bimodule $B$, we have $\H_i(A,B)=0$ for $i>1$ and
 $$\H_0(A,B)\cong \cok(u) \ {\rm and} \ \ \H_1(A,B)\cong \ker(u),$$
 where $u:B\to B$ is given by $u(b)=xb-bx$. If $B=\Phi^n(\F_p[x])$, then
$u:\F_p[x]\to \F_p[x]$ is the multiplication by $(x^{p^n}-x)$ and we
obtain $\H_1(A,\Phi^n(A))=0$ and $\H_0(A,\Phi^n(A))=\psi^n(A)$

{\bf Step 2. The theorem holds if  $A$ is a polynomial algebra}.
Since Hochschild homology commutes with filtered colimits it
suffices to consider the case when $A=\F_p[x_1,\cdots, x_d]$. By
 the K\"unneth theorem
for Hochschild homology (see \cite[Theorem X.7.4]{homology} we have
$\H_\bu(A,\Phi^n(A))=\H_\bu(\F[x],\Phi^n(\F[x]))^{\t d}$ and the
result follows.

{\bf Step 3. The theorem holds for arbitrary $A$}. We use the same
method as  used in the proof of \cite[Theorem 3.5.8]{CH}. First we
choose a simplicial commutative algebra $L_\bu$ such that each $L_n$
is a polynomial algebra, $n\geq 0$ and $\pi_i(L_\bu)=0$ for all
$i>0$, $\pi_0(L_\bu)=A$. Such a resolution exist thanks to
\cite{HA}. Now consider the bisimplicial vector space
$C_\bu(L_\bu,\Phi^n(L_\bu))$.
 The $s$-th horizontal simplicial vector space
is the simplicial vector space $L_\bu^{\t s+1}$. By the
Eilenberg-Zilber-Cartier and K\"unneth theorems it  has zero
homotopy groups in positive dimensions  and $\pi_0(L_\bu^{\t
s+1})=A^{\t s+1}$. On the other hand the $t$-th vertical simplicial
vector space of $C_\bu(L_\bu,\Phi^n(L_\bu)))$ is isomorphic to the
Hochschild complex $C_\bu(L_t,\Phi^n(L_t)))$ which has zero homology
in positive dimensions  by the previous step. Hence both spectral
sequences corresponding to the bisimplicial vector space
$C_\bu(L_\bu,\Phi^n(L_\bu))$ degenerate and we obtain the
isomorphism
$$\H_\bu(A,\Phi^n(A))\cong \pi_\bu(\psi^n(L_\bu)).$$
Now we can use  Lemma \ref{homonulia}  to finish the proof.
\end{proof}
\begin{Co} Let $A$ be a commutative $\F_p$-algebra,  $M$ be an
$A$-module and $n>1$. Then there exist functorial isomorphisms
$$\H_\bu(A,\Phi^n(M))\cong \tor^A_\bu(\psi^n(A),M), \ \ n\geq 0$$
and
$$\H^\bu(A,\Phi^n(M))\cong \ext_A^\bu(\psi^n(A),M), \ \ n\geq 0.$$

In particular, if $A$ is  a commutative algebra over a  field $K$ of
characteristic $p>0$ with ${\sf Card}(K)>p^n$, then
$$\H_\bu(A,\Phi^n(M))=0=\H^\bu(A,\Phi^n(M)).$$
\end{Co}

\begin{proof} Observe that $C_\bu(A,\Phi^n(A))$ is a complex of left
$A$-modules. By  Theorem \ref{respsi} it is a free-resolution of
$\psi^n(A)$ in the category of $A$-modules. Hence it suffices to
note that
$$C_\bu(A,\Phi^n(M))\cong M\t_AC_\bu(A,\Phi^n(A)),$$
$$C^\bu(A,\Phi^n(M))\cong \hom_A(C_\bu(A,\Phi^n(A)),M),$$
where $C^*$ denotes the standard complex for Hochschild cohomology.
The last assertion follows from Lemma \ref{us}.
\end{proof}
\begin{Exm}
It follows for instance that $\H^i(A,\Phi^n(M))=0$, $i>0$, provided
$M$ is an injective $A$-module and $n>1$. In particular
$\H^i(A,\Phi^n(A))=0$ if $A$ is a self-injective algebra.  On the
other hand if $A=\F_p[x_1,\cdots, x_d]$ then $\H^i(A,\Phi^n(A))=0$,
$i\not =d$, $n>1$ and $\H^d(A,\Phi^n(A))=\psi^n(A)$, $n>1$.
\end{Exm}
\section{Application to MacLane cohomology}
We recall the definition of MacLane (co)homology. For  an
associative  ring $R$ we let ${\bf F}(R)$ be the category of
finitely generated free left $R$-modules. Moreover, we let $\sF(R)$
be the category of all covariant functors from the category ${\bf
F}(R)$  to the category of all $R$-modules. The category $\sF(R)$ is
an abelian category with enough projective and injective objects. By
definition \cite{JP} the \emph{MacLane cohomology} of $R$ with
coefficient in a functor $T\in \sF(R)$ is given by
$$\HML^\bu(R,T):=\ext^\bu_{\sF(R)}(I,T),$$
where  $I\in \sF(R)$ is  the inclusion of the category ${\bf F}(R)$
into the category of all left $R$-modules. One defines MacLane
homology in  a dual manner (see \cite[Proposition 3.1]{PW}). For an
$R$-$R$-bimodule $B$, one considers  the functor $B\t_ R(-)$ as an
object of the category $\sF(R)$. For simplicity we write
$\HML_\bu(R,B)$ instead of $\HML_\bu(R,B\t_R(-))$. There is a
binatural transformation
$$\HML_\bu(R,B)\to \H_\bu(R,B)$$ which is an isomorphism in
dimensions $0$ and $1$.

In the rest of this section we consider MacLane (co)homology of
commutative $\F_p$-algebras.

\begin{Le} For any commutative $\F_p$-algebra $A$ one has an
isomorphism
$$\HML_{2i}(A,\Phi^n(A))=\psi^n(A), \ i\geq 0, n>1,$$
and $$\HML_{2i+1}(A,\Phi^n(A))=0, \ i\geq 0, n>1.$$
\end{Le}
\begin{proof} According to \cite[Proposition 4.1]{P}
there exists a functorial spectral sequence
$$E^2_{pq}=\H_p(A,\HML_q(\F_p,B))\Longrightarrow
\HML_{p+q}(A,B).$$ Here $B$ is an $A$-$A$-bimodule. By the
well-known computation of Breen \cite{breen},  B\"okstedt
\cite{marcel} (see also \cite{FLS}) we have
$$\HML_{2i}(\F_p,B)=B$$ and
$$\HML_{2i+1}(\F_p,B)=0.$$
Now we put $B=\psi^n(A)$ and use Theorem \ref{respsi} to get
$E^2_{pq}=0$ for all $p>0$. Hence the spectral sequence degenerates
and the result follows.
\end{proof}

We now consider MacLane cohomology with coefficients in strict
polynomial functors \cite{FS}. Let us recall that the strict
homogeneous polynomial functors of degree $d$ form an abelian
category $\sP_d(A)$ and there exist an exact functor $i:\sP_d(A)\to
\sF(A)$ \cite{FFSS}. For an object $T\in \sP_d(A)$  we write
$\HML_\bu(A,T)$ instead of $\HML_\bu(A,i(T))$. Projective generators
 of the category $\sP_d$ are tensor products of the divided powers,
 while the injective cogenerators are symmetric powers.
 Let us
recall that the $d$-th divided power functor $\Gg^d\in \sF(A)$ and
$d$-th symmetric functors $S^n$ are defined by $$\Gg^d(M)=(M^{\t
d})^{\Sigma_d}, \ \ S^n(M)=(M^{\t d})_{\Sigma_d}.$$ Here tensor
products are taken over $A$, $\Sigma_d$ is the symmetric group on
$d$-letters, which acts on the $d$-th tensor power by permuting of
factors, $M\in {\bf F}(A)$ and $X^G$ (resp. $X_G$) denotes the
module of invariants (resp. coinvariants) of a $G$-module $X$, where
$G$ is a  group.

For a functor $T\in \sF(A)$ we let $\tilde{T}\in \sF(\F_p)$ be the
functor defined by
$$\tilde{T}(V)=T(V\t A)$$
According to \cite[Theorem 4.1]{PW} the groups
$\HML_i(\F_p,\tilde{T})$ have an $A$-$A$-bimodule structure. The
left action comes from the fact that $T$ has values in the category
of left $A$-modules, while the right action  comes from the fact
that $T$ is defined on ${\bf F}(A)$. In particular it uses the
action of $T$ on the maps $l_a:X\to X$, where $a\in A$, $X\in {\bf
F}(A)$ and $l_a$ is the multiplication on $a$. Since
$T(l_a)=l_{a^d}$ if $T$ is a strict homogeneous polynomial functor
of degree $d$ \cite{FS}, the bimodule $\HML_i(\F_p,\tilde{T})$ is of
the form $\Phi^n(M)$ provided $d=p^n$.

\begin{Th} Let $d>1$ be an integer and let $A$ be a commutative
$\F_p$-algebra. Then $\HML_\bu(A,\Gg^d)=0$ if $d$ is not a power of
$p$. If $d=p^n$ and $n>0$, then
$$\HML_i(A,\Gg^d)=0 \ \ {\rm if} \ i\not
=2p^nt, t\geq 0$$ and $$\HML_i(A,\Gg^d)=\psi^n(A)\ {\rm if} \ \
i=2p^nt, t\geq 0.$$ In particular $\HML_\bu(A,\Gg^d)=0$ provided $A$
is an algebra over  a  field $K$ of characteristic $p>0$ with ${\sf
Card}(K)>d$.
\end{Th}
\begin{proof} According to \cite[Theorem 4.1]{PW},\cite{P}
there exists a functorial spectral sequence:
$$E^2_{pq}=\H_p(A,\HML_q(\F_p,\tilde{T}))\Longrightarrow
\HML_{p+q}(A,T).$$ For $T=\Gg_A^n$ one has $\tilde{T}=\Gg^n_{\F_p}\t
A$. Here we used the notation $\Gg_A^n$ in order to emphasize the
dependence  on the ring $A$. By the result of Franjou, Lannes and
Schwartz \cite{FLS}  $\HML_i(\F_p,\tilde{T})$ vanishes unless
$d=p^n$ and $i=2p^nt$, $t\geq 0$. Moreover in these exceptional
cases $\HML_i(\F_p,\tilde{T})$ equals  to $\Phi^n(A)$ (as an
$A$-$A$-module). Hence the spectral sequence together with Theorem
\ref{respsi} gives the result.
\end{proof}

\begin{Co} Let $A$ be  a commutative algebra
over a  field $K$ of characteristic $p>0$ with ${\sf Card}(K)>d>1$.
If $T$ is a strong homogeneous polynomial functor of degree $d$.
Then
$$\HML_\bu(A,T)=0=\HML^\bu(A,T).$$
\end{Co}

\begin{proof} We already proved that the result is true if $T$ is a divided
power. By the well-known vanishing result  \cite{pi} the result is
also true if $T=T_1\t T_2$ with $T_1(0)=0=T_2(0)$. Since any object
of $\sP_d$ has a finite resolution which consists with finite direct
sums of tensor products of divided powers \cite{FS}  the result
follows.
\end{proof}

\end{document}